\documentclass[11pt,letterpaper]{article}

\usepackage{latexsym,multirow}
\usepackage{amssymb,amsmath,bm}
\usepackage{graphicx}
\usepackage[color,all,import,arrow]{xy}
\usepackage{enumerate}
\usepackage{booktabs}
\usepackage{enumitem}

\usepackage{natbib}
\usepackage{xcolor}
\usepackage[pdftex, bookmarksopen=true, bookmarksnumbered=true,
pdfstartview=FitH, breaklinks=true, urlbordercolor={0 1 0}, citebordercolor={0 0 1}, colorlinks=true, citecolor=blue!60!black]{hyperref}
\usepackage{colortbl}
\usepackage{subfigure}
\usepackage{subcaption}

\usepackage{tikz}
\usepackage{inputenc}
\usetikzlibrary{shapes,arrows,trees,fit,positioning,shapes.misc,calc}

\usepackage{dcolumn}
\newcolumntype{.}{D{.}{.}{-1}}
\newcolumntype{d}[1]{D{.}{.}{#1}}

\usepackage{theorem}
\theoremstyle{plain}
\theoremheaderfont{\scshape}
\newtheorem{assumption}{Assumption}
\newtheorem{example}{Example}
\newtheorem{definition}{Definition}
\newtheorem{corollary}{Corollary}

\newtheorem{theorem}{Theorem}

\newtheorem{lemma}{Lemma}
\newenvironment{proof}{\paragraph{Proof:}}{\hfill$\square$}

\newcounter{examplerevisit}
\newenvironment{continuedexample}[2]{%
    \stepcounter{examplerevisit}%
    \renewcommand{\theexample}{\getrefnumber{#1}}%
    \begin{example}[#2, continued]%
}{%
    \end{example}%
    \addtocounter{example}{-1}%
}

\newcommand{\qed}{\hfill \ensuremath{\Box}}
\newcommand{\ind}{\mbox{$\perp\!\!\!\perp$}}

\usepackage{rotating}

\usepackage[compact]{titlesec}

\allowdisplaybreaks

\newcommand\spacingset[1]{\renewcommand{\baselinestretch}%
  {#1}\small\normalsize}

\newcommand{\blind}{0}

\newcommand{\bp}{\bm{p}}
\newcommand{\bq}{\bm{q}}
\newcommand{\bv}{\bm{v}}

\newcommand{\cT}{\mathcal{T}}

\newcommand{\cY}{\mathcal{Y}}

\newcommand{\by}{\bm{y}}

\newcommand{\E}{\mathbb{E}}
\newcommand{\R}{\mathbb{R}}

\newcommand{\bG}{\bm{G}}
\newcommand{\bM}{\bm{M}}

\newcommand{\bw}{\bm{w}}

\newcommand{\cC}{\mathcal{C}}
\newcommand{\cS}{\mathcal{S}}

\usepackage{dsfont}
\usepackage{caption}
\usepackage{tikz-cd}
\usetikzlibrary{arrows.meta}

\newcommand{\lt}{\left}
\newcommand{\rt}{\right}
\newcommand{\indicator}[1]{\mathds{1}{\left\{ #1 \right\}}}

\newcommand{\abs}[1]{\left|#1\right|}

\newcommand{\bYpot}{\bm{Y}^{\ast}}

\newcommand{\cPstar}{\mathcal{P}^{\ast}}

\begin{document} 

\newcommand{\tit}{A General Framework for Identifying the Joint Potential Outcome Distribution through \\ Support Restrictions}
%
%
\spacingset{1.25}

\if0\blind

{\title{\bf\tit}

\author{Benedikt Koch\thanks{Ph.D. candidate, Department of Statistics, Harvard University. 33 Oxford Street, Cambridge MA 02138. Email: \href{mailto:benedikt_koch@g.harvard.edu}{benedikt\_koch@g.harvard.edu} URL: \href{https://benediktjkoch.github.io}{https://benediktjkoch.github.io}} \and Kosuke
  Imai\thanks{Edith and Benjamin Geisinger Professor of Government and of
    Statistics, Harvard University.  1737 Cambridge Street,
    Institute for Quantitative Social Science, Cambridge MA 02138, U.S.A.
    Email: \href{mailto:imai@harvard.edu}{imai@harvard.edu} URL:
    \href{https://imai.fas.harvard.edu}{https://imai.fas.harvard.edu}}
}

\date{
\today
}

\maketitle

}\fi

\if1\blind
\title{\bf \tit}

\maketitle
\fi

\pdfbookmark[1]{Title Page}{Title Page}

\thispagestyle{empty}
\setcounter{page}{0}
         
\begin{abstract}
Many causal estimands depend on the joint distribution of potential outcomes, which is generally not identifiable because only one potential outcome can be observed for each unit. To obtain identification, researchers often impose assumptions such as monotonicity and exclusion restriction that constrain the support of this joint distribution.  For treatments and outcomes that each take finitely many values, we define a similarity matrix whose entries count, for each pair of potential outcome strata, the number of treatments yielding the same outcome.  We then prove that, under unconfoundedness, identification holds if and only if this matrix is invertible, yielding a criterion directly verifiable from the support.  In addition, we derive an identification formula for the joint distribution under these support restrictions and show that an average causal estimand is identifiable if and only if its integrand admits an additive representation on the restricted support.  Finally, we derive sharp testable implications of the support restriction.  Our results unify and generalize a broad range of identification results in the literature.

\end{abstract}

\newpage
\section{Introduction}\label{sec:introdcution}

Many causal quantities of interest depend on the joint, rather than marginal, distribution of potential outcomes.
They include probabilities of causation \citep{robins1989probability, balke1994counterfactual,pearl1999probabilities,tian2000probabilities}, distribution of individual treatment effects \citep{heckman1997making,fan2010sharp}, probabilities of benefit and harm \citep{gadbury2000unit,kallus2022s,mueller2023personalized,ben-michaelPolicyLearningAsymmetric}, and counterfactual utilities \citep{koch2025statistical, koch2026axiomaticfoundationdecisionscounterfactual, imai2026triage}.
However, because we can only observe one potential outcome for every unit, these causal estimands generally remain unidentifiable even when the treatment assignment is unconfounded and the overlap assumption is satisfied.

One common approach to this identification problem is to restrict the support of the joint potential outcome distribution through assumptions such as monotonicity and exclusion restriction.
For binary treatments and outcomes, it is well known that monotonicity identifies the joint potential outcome distribution  \citep{pearl1999probabilities, tian2000probabilities}. For binary outcomes, this identification result generalizes to ordered multi-valued treatments under the assumption of monotone treatment response \citep{manski1997monotone}. \citet{heckman2018unordered} extend monotonicity to unordered choices.  In instrumental variable models, monotonicity and the exclusion restriction together identify the average complier causal effect \citep{imbens1994identification,angrist1995twostage}.  More recently, for ordinal outcomes, \cite{zhang2025identifying} and \cite{hashimoto2026bounds} consider identification under the assumption of monotonic incremental treatment effects.

In this paper, we generalize these identification results by deriving a necessary and sufficient condition for identification of the joint potential outcome distribution through support restrictions. 
Specifically, we first define a similarity matrix whose entry counts the number of treatment values under which a pair of potential outcome strata yields the same outcome.
We then prove that the joint potential outcome distribution is identifiable if and only if this similarity matrix is invertible, and provide an identification formula.
We also show that an average causal quantity, which is a linear functional of the joint potential outcome distribution, can be identified if and only if its integrand admits an additive representation in the potential outcomes on the restricted support \citep[see][for a related result]{koch2025statistical}.
Finally, our identification condition implies two simple bounds on the cardinality of restricted support.
Although these support restrictions cannot be verified directly from the observed data, we derive sharp testable implications that characterize its compatibility with the observed distribution. A violation of these conditions falsifies the support restriction.

Our framework encompasses existing support restrictions as special cases \citep[e.g.,][]{angrist1996identification, manski1997monotone, pearl1999probabilities, tian2000probabilities, blair2012statistical, heckman2018unordered, zhang2025identifying, hashimoto2026bounds}, and generalizes earlier identification conditions, including the sufficient invertibility condition of \citet{pearl1999probabilities} for Markovian models, and the rank condition of \citet{heckman2018unordered} for identifying joint treatment probabilities based on zero--one matrices. In contrast to these works, we derive a necessary and sufficient condition for the identification of the joint potential outcome distribution through the invertibility of a similarity matrix, which is interpretable in terms of Hamming distances between principal strata. We also provide closed-form expressions of identification formulae and obtain sharp testable implications that characterize the compatibility of support restrictions with the observed distribution.

\section{Identification through Support Restrictions}

\subsection{Setup}
Let $T \in \cT := \{0, 1, \ldots, K-1\}$ denote a treatment variable and $Y \in \cY := \{0, 1, \ldots, L-1\}$ be an observed outcome variable, where $K, L \geq 2$. Both treatment and outcome variables take finitely many values, which may or may not be ordered. We also define the vector of potential outcomes under all treatment values as $\bYpot := (Y(0), \ldots, Y(K-1)) \in \cY^{K}$. The full data vector $(T, \bYpot)$ follows an unknown super-population law $P^*$. We link the full data vector to the observed data vector through the standard consistency assumption.
\begin{assumption}[Consistency]\label{ass:consistency}
    The observed outcome variable satisfies $Y = Y(T)$.
\end{assumption}

Under Assumption~\ref{ass:consistency}, the observed data vector $(T, Y)$ have law $P = P^* \circ \cC^{-1}$, where $\cC(t, \by) = (t, y_t)$ is a coarsening map. We impose the following standard identification assumption on $P^*$ \citep{rosenbaum1983central}. For notational simplicity, we do not condition on a set of observed pre-treatment covariates.  However, all of our results hold, conditional on the covariates.
\begin{assumption}[Strong Ignorability]\label{ass:ignorability}
\mbox{}\par
\begin{enumerate}[label=(\alph*)]
\item \label{ass:unconfoundedness} Unconfoundedness: $T \ind \bYpot$.
\item \label{ass:overlap} Overlap: $P(T = t) > 0$ for every $t \in \cT$.
\end{enumerate}
\end{assumption}

Under Assumptions~\ref{ass:consistency}~and~\ref{ass:ignorability} we can write the marginal potential outcome probabilities in terms of the observed data law, i.e.,
\[
    P^*(Y(t) = y) = P(Y = y \mid T = t),
\]
for all $t \in \cT$ and $y \in \cY$. However, our goal is to derive the conditions under which we can identify the full joint law $P^*$ from $P$. To this end, define the set of joint potential outcome distributions that satisfy Assumption~\ref{ass:ignorability},
\[
    \cPstar := \{P^*: P^* \text{ satisfies Assumption~\ref{ass:ignorability}}\}.
\]
In addition, we use the following standard notion of identifiability \citep{pearl1999probabilities}.
\begin{definition}[Identifiability]\label{def:identifiability}
    Let $\cPstar_0 \subseteq \cPstar$. An estimand $\theta$ defined on $\cPstar$ is \emph{identifiable} over $\cPstar_0$ if, for every $P_1^*, P_2^* \in \cPstar_0$ that induce the same observed law, $P_1 = P_2$, we have $\theta(P_1^*) = \theta(P_2^*)$. 
\end{definition}

Thus, identification requires that two full data laws that generate the same observed law, agree on the estimand $\theta$. We first study identification of the full-data law itself, $\theta(P^\ast) = P^\ast$, and then consider expectations of functions under this law. Because only one potential outcome can be observed for each unit, Assumption~\ref{ass:ignorability} generally does not suffice for identification \citep{holland1986}. 

We consider a set of identification assumptions that restrict the the support of $\bYpot$. As discussed in Section~\ref{sec:introdcution}, this rerepsents one of the most common identification approaches in the literature.  
\begin{assumption}[Restricted Support]\label{ass:support}
    The potential outcome vector satisfies $\bYpot \in \cS$ almost surely where $\cS \subseteq \cY^{K}$ is a nonempty set.
\end{assumption}

One prominent example of Assumption~\ref{ass:support} is monotonicity \citep[e.g.,][]{manski1997monotone,pearl1999probabilities, zhang2025identifying, hashimoto2026bounds}. For example, for binary treatments and outcomes, monotonicity requires that $Y(0) \leq Y(1)$, restricting the support of $\bYpot = (Y(0), Y(1)) \in \cS = \{(0,0),(0,1),(1,1)\}$. Throughout, we apply Definition~\ref{def:identifiability} to the class of full data laws that satisfy a given support restriction $\cS$, 
\[
    \cPstar_{\cS}:= \{P^* \in \cPstar: P^* \text{ satisfies Assumption~\ref{ass:support} with } \cS\}.
\]

\subsection{Necessary and Sufficient Condition for Identification}

We demonstrate that simple matrix inversion fully characterizes when a given support restriction identifies the full data law $P^*$. 
For the restricted support $\cS$ defined in Assumption~\ref{ass:support}, we construct a $|\cS|\times|\cS|$ integer-valued ``similarity'' matrix $\bG_{\cS}$ whose $(\by,\by^\prime)$ entry counts the treatment values for which the two potential outcome vectors agree:
\begin{equation}\label{eq:gram_similarity}
    (\bG_{\cS})_{\by,\by^\prime} :=  \sum_{t \in \cT} \indicator{y_t = y_t^\prime} \ = \ K - d_H(\by,\by^\prime),
\end{equation}
for $\by,\by^\prime \in \cS $ where $d_H(\by,\by^\prime) := \sum_{t \in \cT} \indicator{y_t \neq y_t^\prime}$ is the Hamming distance.

The following theorem shows that the inversion of this similarity matrix is a necessary and sufficient condition for the identification of joint potential outcome distribution.  The inverted matrix also provides a closed-form identification formula for the expected value of any function under $P^*$.
\begin{theorem}[Identification of the Joint Potential Outcome Distribution]\label{thm:joint_identification}
    The law $P^*$ is identifiable over $\cPstar_{\cS}$ if and only if $\bG_{\cS}$ is invertible. In that case, for every $P^* \in \cPstar_{\cS}$ and $\by \in \cS$,
    \begin{equation}\label{eq:recovery}
        P^*(\bYpot = \by) = \sum_{\by^\prime \in \cS} (\bG_{\cS}^{-1})_{\by,\by^\prime} \sum_{t \in \cT} P(Y = y_t^\prime \mid T = t).
    \end{equation}
    In particular, for any function $h:\cT\times\cY^{K}\to\R$,
    \begin{align*}
    \E_{P^*}[h(T,\bYpot)] &= \sum_{\by,\by^\prime \in \cS} (\bG_{\cS}^{-1})_{\by,\by^\prime} \sum_{t \in \cT} h(t,\by)P(T=t) \sum_{k \in \cT} P(Y = y_k^\prime \mid T = k).
\end{align*} 
\end{theorem}
The proof of Theorem~\ref{thm:joint_identification} is given in Appendix~\ref{app:joint_identification}.

Equation~\eqref{eq:recovery} identifies the principal stratum proportions by inverting the following equation,
\[
    \sum_{t \in \cT} P(Y = y_t^\prime \mid T = t)
    = \sum_{\by \in \cS} (\bG_{\cS})_{\by^\prime,\by} P^*(\bYpot = \by),
\]
for $\by^\prime \in \cS$. In other words, each observable sum of marginal probabilities is a weighted sum of stratum proportions, with weights given by the Hamming similarities. If $\bG_{\cS}$ is not invertible, its columns are linearly dependent. Redistributing probability mass among strata can therefore leave all observable sums unchanged, so that multiple full-data laws can induce the same observed law. Invertibility of $\bG_{\cS}$ ensures that this mapping is unique. Below, we illustrate the application of this theorem through a simple example of monotonicity with binary treatment and outcome \citep{pearl1999probabilities}. 

\begin{example}[Monotonicity with binary treatment and outcome]\label{ex:binary}
    Let $\cT = \cY = \{0, 1\}$. Monotonicity imposes $Y(0) \leq Y(1)$ almost surely. Hence the support of $\bYpot = (Y(0), Y(1))$ is given by $\cS = \{(0,0),(0,1),(1,1)\}$, excluding $(1, 0)$. The corresponding similarity matrix $\bG_{\cS}$ in Equation~\eqref{eq:gram_similarity} is invertible, with
    \[
        \bG_{\cS}
            = \begin{pmatrix}
                2 & 1 & 0\\
                1 & 2 & 1\\
                0 & 1 & 2
              \end{pmatrix}, \qquad \bG_{\cS}^{-1}
        = \frac{1}{4}
          \begin{pmatrix}
              3 & -2 & 1\\
              -2 & 4 & -2\\
              1 & -2 & 3
          \end{pmatrix}.
\]
    Theorem~\ref{thm:joint_identification} implies that $P^*$ is identifiable. Write $q_t := P(Y=1\mid T=t)$.  Equation~\eqref{eq:recovery} yields
    \begin{align*}
        P^*(\bYpot=(0,0))
        &= \frac{1}{4}\lt[3(2-q_0-q_1)-2(1-q_0+q_1)+(q_0+q_1)\rt] = 1-q_1 \\
        &= P(Y=0\mid T=1),\\
        P^*(\bYpot=(0,1))
        &= \frac{1}{4}\lt[-2(2-q_0-q_1)+4(1-q_0+q_1)-2(q_0+q_1)\rt] = q_1-q_0 \\
        &= P(Y=1\mid T=1)-P(Y=1\mid T=0),\\
        P^*(\bYpot=(1,1))
        &= \frac{1}{4}\lt[(2-q_0-q_1)-2(1-q_0+q_1)+3(q_0+q_1)\rt] = q_0 \\
        &= P(Y=1\mid T=0).
    \end{align*}
\end{example}

The support restriction in Example~\ref{ex:binary} extends to finitely many treatments and outcomes under monotone treatment response \citep{manski1997monotone}. However, as the next example shows, this restriction identifies the joint distribution only when the outcome is binary.

\begin{example}[Monotone treatment response]\label{ex:mtr}
    Monotone treatment response requires $Y(s)\leq Y(t)$ almost surely whenever $s\leq t$. First, let  $\cT = \{0, 1\}$ and $\cY = \{0, 1, 2\}$. Then 
    \[
        \cS=\{(0,0),(0,1),(0,2),(1,1),(1,2),(2,2)\}
    \]
    and
    \[
        \bG_{\cS}
        = \begin{pmatrix}
            2 & 1 & 1 & 0 & 0 & 0\\
            1 & 2 & 1 & 1 & 0 & 0\\
            1 & 1 & 2 & 0 & 1 & 1\\
            0 & 1 & 0 & 2 & 1 & 0\\
            0 & 0 & 1 & 1 & 2 & 1\\
            0 & 0 & 1 & 0 & 1 & 2
          \end{pmatrix}.
    \]
    Its columns are linearly dependent, since
    \[
        \bG_{\cS}[:,(0,1)] + \bG_{\cS}[:,(1,2)] = \bG_{\cS}[:,(0,2)] + \bG_{\cS}[:,(1,1)].
    \]
    Thus, by Theorem~\ref{thm:joint_identification}, $P^*$ is not identifiable over $\cPstar_{\cS}$. The same conclusion holds for every $L = |\cY| \geq 3$, since the support must contain $\cS$.

    Assume now binary outcomes $\cY = \{0, 1\}$ but a general treatment with $K = |\cT| \geq2$. Then 
    \[
        \cS = \{\by \in \{0,1\}^K : y_0 \leq \cdots \leq y_{K-1}\} = \{(0,\ldots,0),\ (0,\ldots,0,1),\ \ldots,\ (0,1,\ldots,1),\ (1,\ldots,1)\}.
    \]
    Each stratum $\by \in \cS$ is uniquely determined by its number of ones, $i:=\sum_{t \in \cT} y_t$. Strata indexed by $i$ and $j$ differ under exactly $|i-j|$ treatments, so Equation~\eqref{eq:gram_similarity} gives
    \[
    (\bG_{\cS})_{i,j}=K-|i-j|,
    \]
    for $i, j = 0, \ldots, K$. Since columns of $\bG_{\cS}$ are linearly independent, identification holds under monotone treatment response if and only if $L=2$ for any $K\geq2$. In this case, the stratum probabilities can be read directly from the support. Let $\by^{(i)}$ denote the stratum with exactly $i$ ones. Then
    \[
        P^*(\bYpot=\by^{(i)}) =
    \begin{cases}
        P(Y=0\mid T=K-1),                          & i=0,\\
        P(Y=1\mid T=K-i)-P(Y=1\mid T=K-i-1),      & 1\leq i\leq K-1,\\
        P(Y=1\mid T=0),                            & i=K.
    \end{cases}
    \]
\end{example}

Example~\ref{ex:mtr} shows that monotonicity alone does not identify the joint distribution for non-binary outcomes. As the next example shows, identification holds if we additionally assume that the treatment can increase each unit's outcome by at most one level \citep{zhang2025identifying,hashimoto2026bounds}.

\begin{example}[Monotonic incremental treatment effects]\label{ex:mite}
    Monotonic incremental treatment effects require that $0\leq Y(t)-Y(s)\leq1$ whenever $s\leq t$. That is, as the treatment increases, each unit's outcome can increase at most once, by one level. This gives the support
    \[
        \cS = \{\by\in\cY^{K}:y_0\leq\cdots\leq y_{K-1},\ y_{K-1}-y_0\leq1\}.
    \]
    Compared with Example~\ref{ex:mtr}, when $K=2$ and $L=3$, this excludes stratum $(0,2)$ from $\cS$. For general $K,L\geq2$, each stratum $\by\in\cS$ is again uniquely determined by the sum of its entries, $i:=\sum_{t\in\cT}y_t$. Strata indexed by $i$ and $j$ have $K-|i-j|$ matching entries when $|i-j|<K$, and none otherwise. The corresponding similarity matrix is given by
     \[
        (\bG_{\cS})_{i, j} = \max\{K-|i-j|,0\},
    \]
    for $i, j = 0, \ldots, K(L-1)$. Since its columns are linearly independent, $\bG_{\cS}$ is invertible, and Theorem~\ref{thm:joint_identification} implies that $P^*$ is identifiable.

    The stratum probabilities can also be read directly from the support. Each stratum $\by \in \cS$ either remains at $y_0 \in \cY$ or increases by one at a unique treatment. Hence,
    \begin{equation}\label{eq:mite_stratum_probabilities}
        P^*(\bYpot=\by)=
        \begin{cases}
            P(Y\geq y_0\mid T=0)-P(Y>y_0\mid T=K-1), & y_{K-1}=y_0,\\
            P(Y>y_0\mid T=t)-P(Y>y_0\mid T=t-1), & y_t>y_{t-1}.
        \end{cases}
    \end{equation}
\end{example}

To accommodate settings in which the outcome categories are labels rather than ordered, \citet{heckman2018unordered} propose unordered monotonicity. The next example shows that this condition is generally not strong enough to guarantee identification. 

\begin{example}[Unordered monotonicity]\label{ex:unordered}
    Unordered monotonicity imposes that, for every $s,t\in\cT$, $a\in\cY$, and $\by\in\cS$, we have either
    $$
    \indicator{y_s=a}\leq\indicator{y_t=a}
        \quad \text{ or } \quad
    \indicator{y_s=a}\geq\indicator{y_t=a}.
    $$
    This restriction does not generally identify the full data law $P^*$. For example, when $K=2$ and $L=4$, the support $\cS=\{(0,0),(1,1),(2,2),(3,3),(0,2),(0,3),(1,2),(1,3)\}$
    satisfies unordered monotonicity. Its similarity matrix is
    \[
        \bG_{\cS}
        = \begin{pmatrix}
            2 & 0 & 0 & 0 & 1 & 1 & 0 & 0\\
            0 & 2 & 0 & 0 & 0 & 0 & 1 & 1\\
            0 & 0 & 2 & 0 & 1 & 0 & 1 & 0\\
            0 & 0 & 0 & 2 & 0 & 1 & 0 & 1\\
            1 & 0 & 1 & 0 & 2 & 1 & 1 & 0\\
            1 & 0 & 0 & 1 & 1 & 2 & 0 & 1\\
            0 & 1 & 1 & 0 & 1 & 0 & 2 & 1\\
            0 & 1 & 0 & 1 & 0 & 1 & 1 & 2
          \end{pmatrix}.
    \]
    Its columns are linearly dependent, since
    \[
        \bG_{\cS}[:,(0,2)] + \bG_{\cS}[:,(1,3)]
        = \bG_{\cS}[:,(0,3)] + \bG_{\cS}[:,(1,2)].
    \]
    Thus, by Theorem~\ref{thm:joint_identification}, $P^*$ is not identifiable over $\cPstar_{\cS}$.
\end{example}

Although all preceding examples impose a form of monotonicity, this is not the only type of support restriction that can be used for identification, as illustrated by Example~\ref{ex:noisy_or} in Appendix~\ref{app:examples}.


\subsection{Simple Bounds on the Support Size}

Although Theorem~\ref{thm:joint_identification} provides a necessary and sufficient condition for identification, the similarity matrix $\bG_{\cS}$ can be high-dimensional when the restricted support $\cS$ is large. The following two simple bounds on the support size provide necessary conditions for identification that are easy to check.
\begin{corollary}[Cardinality Constraints on the Support]\label{cor:cardinality}
    For each $t \in \cT$, let $\cS_t := \{y_t:\by \in \cS\}$. Identifiability of $P^*$ over $\cPstar_{\cS}$ requires both
    \begin{enumerate}[label=(\alph*)]
        \item \label{cond:size_refined} $|\cS| \leq 1 + \sum_{t \in \cT} (|\cS_t|-1)$,
        \item \label{cond:size_simple} $|\cS| \leq 1 + K(L-1)$.
    \end{enumerate}
    These conditions are necessary but not sufficient for identification.
\end{corollary}
The proof of Corollary~\ref{cor:cardinality} is given in Appendix~\ref{app:cardinality}. Condition~\ref{cond:size_simple} follows from Condition~\ref{cond:size_refined} but is easier to verify, as it depends only on $|\cS|$, $K$, and $L$ and not on the individual support of each potential outcome. Failure of either condition implies that $P^*$ is not identified over $\cPstar_{\cS}$. A related bound to Condition~\ref{cond:size_simple} under an instrumental variable model is given by \citet[Theorem~T-7]{heckman2018unordered}.

We illustrate the simplicity of the cardinality bound by revisiting Examples~\ref{ex:mtr}~and~\ref{ex:unordered} where we establish the lack of identification while avoiding the construction of the similarity matrix.

\begin{continuedexample}{ex:mtr}{Monotone treatment response}
    Monotone treatment response restricts the support to $\cS = \{\by\in\cY^{K}:y_0\leq\cdots\leq y_{K-1}\}$. For any $K\geq2$ and $L\geq3$,
    \[
        |\cS|=\binom{K+L-1}{K}>1+K(L-1).
    \]
    Thus, Condition~\ref{cond:size_simple} of Corollary~\ref{cor:cardinality} fails, and $P^*$ is not identifiable over $\cPstar_{\cS}$.
\end{continuedexample}

\begin{continuedexample}{ex:unordered}{Unordered monotonicity}
    Recall that for $K=2$ and $L=4$, unordered monotonicity permits the support $\cS = \{(0,0),(1,1),(2,2),(3,3),(0,2),(0,3),(1,2),(1,3)\}$. Since
    \[
        |\cS|=8>7=1+K(L-1).
    \]
    Condition~\ref{cond:size_simple} of Corollary~\ref{cor:cardinality} is violated, and $P^*$ is not identifiable over $\cPstar_{\cS}$.
\end{continuedexample}

\subsection{Identification of Expectations }\label{subsec:expectations}

Theorem~\ref{thm:joint_identification} presents a necessary and sufficient conditio for identifying the joint potential outcome distribution.  However, many causal estimands are expectations of functions of potential outcomes, which are linear functionals of the joint potential outcome distribution.  Some of these estimands are identifiable even when the joint potential outcome distribution may not.
Building on the result of \citet{koch2025statistical}, the following theorem presents a necessary and sufficient condition for identifying expectations under the restricted support and provides a closed-form identification formula.
\begin{theorem}[Identification of Expectations]\label{thm:expectation_identification}
    Let $h:\cT\times\cY^{K}\to\R$. Then $\E_{P^*}[h(T,\bYpot)]$ is identifiable over $\cPstar_{\cS}$ if and only if there exist functions $h_k:\cT\times\cY\to\R$, $k\in\cT$, such that
    \begin{equation*}
        h(t,\by)=\sum_{k\in\cT}h_k(t,y_k),
    \end{equation*}
    for every $t\in\cT$ and $\by\in\cS$. In that case, for every $P^*\in\cPstar_{\cS}$,
    \begin{equation*}
        \E_{P^*}[h(T,\bYpot)]
        =\sum_{t\in\cT}P(T=t)\sum_{k\in\cT}\sum_{a\in\cY}
        h_k(t,a)P(Y=a\mid T=k).
    \end{equation*}
\end{theorem}

The proof of Theorem~\ref{thm:expectation_identification} follows from \citet[Corollary~1]{koch2025statistical} applied to the restricted support $\cS$. Additivity is required only on $\cS$. Thus, a function that is not additive on $\cY^{K}$ may still have an identifiable expectation under a support restriction, even when $\bG_{\cS}$ is not invertible. When $\bG_{\cS}$ is invertible, every function admits an additive representation on $\cS$, and Theorem~\ref{thm:joint_identification} identifies all such expectations. As an illustration, we apply Theorem~\ref{thm:expectation_identification} to unordered monotonicity \citep{heckman2018unordered}.

\begin{continuedexample}{ex:unordered}{Unordered monotonicity}
    As we saw in Example~\ref{ex:unordered}, unordered monotonicity does not generally identify $P^*$. 
    However, one may be interested in particular stratum probabilities, including the proportion of ``always-takers,'' given here as $P^*(\bYpot=(a,\ldots,a))$ for $a\in\cY$ \citep{heckman2018unordered}. Define
    \[
        h(t,\by):=\indicator{\by=(a,\ldots,a)}.
    \]
    Under unordered monotonicity, there exists $t_a\in\cT$ such that for every $\by \in \cS$,
\begin{equation}\label{eq:unordered_additivity}
        h(t,\by)
        =\indicator{y_{t_a}=a}
        =\sum_{k\in\cT}\indicator{k=t_a}\indicator{y_k=a}
        =\sum_{k\in\cT}h_k(t,y_k),
    \end{equation}
    where $h_k(t,b):=\indicator{k=t_a}\indicator{b=a}$ for $t, k\in\cT$ and $b\in\cY$. To illustrate, recall that for $K = 2$ and $L = 4$, the support in Example~\ref{ex:unordered} is given by
    \[
        \cS=\{(0,0),(1,1),(2,2),(3,3),(0,2),(0,3),(1,2),(1,3)\}.
    \]
    Here, $t_a=1$ for $a\in\{0,1\}$ and $t_a=0$ for $a\in\{2,3\}$.
    By Equation~\eqref{eq:unordered_additivity}, $h$ has an additive representation on $\cS$ and Theorem~\ref{thm:expectation_identification} yields
        \begin{align*}
        P^*(\bYpot=(a,\ldots,a))
        &= \E_{P^*}[h(T,\bYpot)]\\
        &= \sum_{t\in\cT}P(T=t)\sum_{k\in\cT}\sum_{b\in\cY}
            h_k(t,b)P(Y=b\mid T=k)\\
        &= \sum_{t\in\cT}P(T=t)\sum_{k\in\cT}\sum_{b\in\cY}
            \indicator{k=t_a}\indicator{b=a}P(Y=b\mid T=k)\\
        &= \sum_{t\in\cT}P(T=t)P(Y=a\mid T=t_a)\\
        &= P(Y=a\mid T=t_a)\\
        &= \min_{t\in\cT}P(Y=a\mid T=t),
    \end{align*}
    where the last equality follows from $\indicator{y_{t_a}=a}\leq\indicator{y_t=a}$ for every $t\in\cT$ and $\by\in\cS$ under unordered monotonicity. This recovers the result established in \citet[Corollary~C-2]{heckman2018unordered}.
\end{continuedexample}

\section{Sharp Testable Implications}

The support restriction in Assumption~\ref{ass:support} is central to our approach, but we cannot directly verify it from the observed data. In this section, we derive its sharp testable implications, which researchers can use in practice to probe the assumption. These conditions are sharp in the sense that (1) if they are violated, no joint potential outcome distribution supported on $\cS$ can be compatible with the observed data distribution, and (2) if they are met, at least one such distribution exists. 

Suppose that $\bG_{\cS}$ is invertible.  Then, define a candidate law
\begin{equation}\label{eq:candidate_probabilities}
    \pi_{\by}(P) := \sum_{\by^\prime \in \cS} (\bG_{\cS}^{-1})_{\by,\by^\prime} \sum_{t \in \cT} P(Y = y_t^\prime \mid T = t),
\end{equation}
for $\by \in \cS$. Under Assumptions~\ref{ass:consistency}--\ref{ass:support}, Theorem~\ref{thm:joint_identification} implies that such a law gives the true stratum probabilities, $\pi_{\by}(P) = P^*(\bYpot = \by)$ for $\by \in \cS$. When Assumption~\ref{ass:support} fails, however, this need not be the case.  The next theorem derives the sharp testable implications. 

\begin{theorem}[Sharp Testable Implications]\label{thm:testable_implications}
    Assume that $\bG_{\cS}$ is invertible.
    \begin{enumerate}[label=(\alph*)]
        \item Testable Implications:  Suppose that $P^*$ satisfies Assumptions~\ref{ass:consistency}  and~\ref{ass:ignorability}. If Assumption~\ref{ass:support} also holds for $P^*$, then its observed law $P=P^*\circ\cC^{-1}$ obeys,
        \begin{align}\label{eq:testable_conditions}
            \pi_{\by}(P) &\geq 0, \qquad
            \sum_{\by \in \cS: y_t = a} \pi_{\by}(P)
            = P(Y = a \mid T = t),
        \end{align}
        for every $\by \in \cS$, $t \in \cT$, and $a \in \cY$, where $\pi_{\by}(P)$ is defined in Equation~\eqref{eq:candidate_probabilities}. Failure to satisfy Equation~\eqref{eq:testable_conditions} implies a violation of Assumption~\ref{ass:support}.

        \item Sharpness: Let $P$ be an observed law satisfying overlap (Assumption~\ref{ass:ignorability}\ref{ass:overlap}). If Equation~\eqref{eq:testable_conditions} holds for $P$, then there exists $P^* \in \cPstar_{\cS}$ such that $P=P^*\circ\cC^{-1}$.
    \end{enumerate}
\end{theorem}
The proof of Theorem~\ref{thm:testable_implications} is given in Appendix~\ref{app:testable_implications}. The conditions in Theorem~\ref{thm:testable_implications} require that the stratum candidate probabilities are nonnegative and sum to the observed marginals. If both conditions are satisfied, the probabilities sum to one.  As an illustration, we apply this theorem to monotonic incremental treatment effects \citep{zhang2025identifying,hashimoto2026bounds}.

\begin{continuedexample}{ex:mite}{Monotonic incremental treatment effects}
    When $K=2$ and $L=3$, the support under monotonic incremental treatment effects is given by
    \[
    \cS = \{(0,0),(0,1),(1,1),(1,2),(2,2)\}.
    \]
    To verify the testable implications of Assumption~\ref{ass:support} on the support, we derive the candidate stratum probabilities in Equation~\eqref{eq:candidate_probabilities} from Equation~\eqref{eq:mite_stratum_probabilities}
    \begin{align*}
        \pi_{(0,0)}(P) &= P(Y\geq 0\mid T=0)-P(Y>0\mid T=1) = P(Y=0\mid T=1),\\
        \pi_{(0,1)}(P) &= P(Y>0\mid T=1)-P(Y>0\mid T=0) = P(Y=0\mid T=0)-P(Y=0\mid T=1),\\
        \pi_{(1,1)}(P) &= P(Y\geq 1\mid T=0)-P(Y>1\mid T=1) = 1-P(Y=0\mid T=0)-P(Y=2\mid T=1),\\
        \pi_{(1,2)}(P) &= P(Y>1\mid T=1)-P(Y>1\mid T=0) = P(Y=2\mid T=1)-P(Y=2\mid T=0),\\
        \pi_{(2,2)}(P) &= P(Y\geq 2\mid T=0)-P(Y>2\mid T=1) = P(Y=2\mid T=0).
    \end{align*}
    These candidates are consistent with the observed marginals in  Equation~\eqref{eq:testable_conditions}. Thus, the testable implications of  Theorem~\ref{thm:testable_implications} reduce to nonnegativity,
    \begin{align*}
        P(Y=0\mid T=0) \geq P(Y=0\mid T=1), \quad 
        P(Y=2\mid T=1) &\geq P(Y=2\mid T=0), \quad 
        \\
        P(Y=0\mid T=0)+P(Y=2\mid T=1) &\leq1.
    \end{align*}
\end{continuedexample}

\section{Conclusion}

Restricting the support of the joint distribution of potential outcomes is a common identification strategy in causal inference. Under consistency, unconfoundedness, and overlap, we show that this distribution is identified if and only if a matrix determined solely by the restricted support is invertible. The matrix measures similarity between pairs of potential outcome strata using the Hamming distance. When it is invertible, the joint distribution and its linear functionals admit closed-form expressions in terms of observed probabilities. This criterion allows researchers to determine whether support restrictions motivated by domain knowledge suffice for identification.

Moreover, researchers are often interested in expectations of functions of potential outcomes rather than the full joint distribution. We show that such an expectation is identifiable if and only if its integrand admits an additive representation on the restricted support. This condition can hold even when the joint distribution is not identifiable.

Finally, we derive sharp testable implications of the support restriction. The observed data rule out the support restriction when the candidate probabilities obtained by inverting the similarity matrix are negative or inconsistent with the observed outcome distribution.


\newpage
\bibliography{ref}

\newpage
\appendix
\renewcommand\thefigure{\thesection.\arabic{figure}}    
\setcounter{figure}{0}
\numberwithin{figure}{section}
\renewcommand\thetable{\thesection.\arabic{table}}    
\setcounter{table}{0} 
\numberwithin{table}{section}
\renewcommand\theequation{\thesection.\arabic{equation}}    
\setcounter{equation}{0}
\numberwithin{equation}{section}
\renewcommand\theproposition{\thesection.\arabic{proposition}}    
\setcounter{proposition}{0}
\numberwithin{proposition}{section}
\numberwithin{assumption}{section}

\begin{center}
 \huge Supplementary Appendix 
\end{center}

\section{Further Examples of Support Restrictions}\label{app:examples}

\begin{example}[Noisy-OR mechanism; \citealp{pearl1999probabilities}]\label{ex:noisy_or}
    Let $\cT=\{(0,0),(1,0),(0,1),(1,1)\}$ and $\cY=\{0,1\}$. For binary auxiliary variables $U_0,U_1,U_2$, consider the mechanism
    \[
        Y(t)=U_0\vee(t_1\wedge U_1)\vee(t_2\wedge U_2),
    \]
    for $t=(t_1,t_2)\in\cT$, where $\vee$ and $\wedge$ denote logical ``or'' and ``and''. 
    Writing $\bYpot=(Y(0,0),Y(1,0),Y(0,1),Y(1,1))$, this restricts the potential outcomes to
    \[
        \cS=\{(0,0,0,0),(0,0,1,1),(0,1,0,1),(0,1,1,1),(1,1,1,1)\}.
    \]
    The similarity matrix and its inverse are
    \[
        \bG_{\cS}
        = \begin{pmatrix}
            4 & 2 & 2 & 1 & 0\\
            2 & 4 & 2 & 3 & 2\\
            2 & 2 & 4 & 3 & 2\\
            1 & 3 & 3 & 4 & 3\\
            0 & 2 & 2 & 3 & 4
          \end{pmatrix}, \qquad
        \bG_{\cS}^{-1}
        = \frac{1}{8}
          \begin{pmatrix}
               5 & -4 & -4 &  4 &  1\\
              -4 &  8 &  4 & -8 &  0\\
              -4 &  4 &  8 & -8 &  0\\
               4 & -8 & -8 & 16 & -4\\
               1 &  0 &  0 & -4 &  5
          \end{pmatrix}.
    \]
    Hence Theorem~\ref{thm:joint_identification} implies that $P^*$ is identifiable, without requiring independence among $U_0,U_1,U_2$.

        To recover the stratum probabilities, recall Equation~\eqref{eq:recovery},
    \[
        P^*(\bYpot=\by)
        =\sum_{\by^\prime\in\cS}
        (\bG_{\cS}^{-1})_{\by,\by^\prime}
        \sum_{t\in\cT}P(Y=y_t^\prime\mid T=t).
    \]
    Write $q_{ij}:=P(Y=1\mid T=(i,j))$ for $i,j\in\{0,1\}$, and denote the inner sum by
    \[
        r_{\by^\prime}:=\sum_{t\in\cT}P(Y=y_t^\prime\mid T=t),
    \]
    for $\by^\prime\in\cS$. The five inner sums are
    \begin{align*}
        r_{(0,0,0,0)}
        &= (1-q_{00})+(1-q_{10})+(1-q_{01})+(1-q_{11})\\
        &= 4-q_{00}-q_{10}-q_{01}-q_{11},\\
        r_{(0,0,1,1)}
        &= (1-q_{00})+(1-q_{10})+q_{01}+q_{11}\\
        &= 2-q_{00}-q_{10}+q_{01}+q_{11},\\
        r_{(0,1,0,1)}
        &= (1-q_{00})+q_{10}+(1-q_{01})+q_{11}\\
        &= 2-q_{00}+q_{10}-q_{01}+q_{11},\\
        r_{(0,1,1,1)}
        &= (1-q_{00})+q_{10}+q_{01}+q_{11}\\
        &= 1-q_{00}+q_{10}+q_{01}+q_{11},\\
        r_{(1,1,1,1)}
        &= q_{00}+q_{10}+q_{01}+q_{11}.
    \end{align*}
    The resulting stratum probabilities are given by weighing the marginals by the rows of $\bG_{\cS}^{-1}$
    \begin{align*}
        P^*(\bYpot=(0,0,0,0))
        &= \frac{1}{8}\lt[5r_{(0,0,0,0)}
            -4r_{(0,0,1,1)}-4r_{(0,1,0,1)}
            +4r_{(0,1,1,1)}+r_{(1,1,1,1)}\rt]\\
        &= 1-q_{11},\\
        P^*(\bYpot=(0,0,1,1))
        &= \frac{1}{8}\lt[-4r_{(0,0,0,0)}
            +8r_{(0,0,1,1)}+4r_{(0,1,0,1)}
            -8r_{(0,1,1,1)}\rt]\\
        &= q_{11}-q_{10},\\
        P^*(\bYpot=(0,1,0,1))
        &= \frac{1}{8}\lt[-4r_{(0,0,0,0)}
            +4r_{(0,0,1,1)}+8r_{(0,1,0,1)}
            -8r_{(0,1,1,1)}\rt]\\
        &= q_{11}-q_{01},\\
        P^*(\bYpot=(0,1,1,1))
        &= \frac{1}{8}\lt[4r_{(0,0,0,0)}
            -8r_{(0,0,1,1)}-8r_{(0,1,0,1)}
            +16r_{(0,1,1,1)}-4r_{(1,1,1,1)}\rt]\\
        &= q_{10}+q_{01}-q_{11}-q_{00} \\
        P^*(\bYpot=(1,1,1,1))
        &= \frac{1}{8}\lt[r_{(0,0,0,0)}
            -4r_{(0,1,1,1)}+5r_{(1,1,1,1)}\rt]\\
        &= q_{00}.
    \end{align*}
\end{example}

\section{Proofs}\label{app:proofs}

\subsection{Proof of Theorem~\ref{thm:joint_identification}}\label{app:joint_identification}

Let $\Delta(\cS) = \{\bp \in [0, 1]^{\cS}: \sum_{\by \in \cS}p_{\by} = 1\}$ be the probability simplex of $\bYpot$ supported on $\cS$.

First, suppose that the joint law $P^*$ is identifiable over $\cPstar_{\cS}$. Consider the marginalization matrix $\bM_{\cS}$ defined in Lemma~\ref{lem:marginalization}. Let $\bp,\widetilde{\bp}\in\Delta(\cS)$ satisfy
\begin{equation}\label{eq:identical_marginals}
    \bM_{\cS}\bp=\bM_{\cS}\widetilde{\bp}.
\end{equation}
We will establish that $\bp = \widetilde{\bp}$.
Define two full data laws by
\begin{align*}
    P^*(T=t,\bYpot=\by) &= \frac{1}{K} p_{\by}, \\
    \widetilde P^*(T=t,\bYpot=\by)
    &= \frac{1}{K} \widetilde p_{\by},
\end{align*}
for $t\in\cT$ and $\by\in\cS$. Under both laws, $T$ is uniformly distributed and independent of $\bYpot$, so Assumption~\ref{ass:ignorability} holds. Trivially, $P^*(\bYpot \in \cS) =   \widetilde P^*(\bYpot \in \cS) = 1$, so Assumption~\ref{ass:support} holds as well. Thus, $P^*,\widetilde P^*\in\cPstar_{\cS}$.

Let $P=P^*\circ\cC^{-1}$ and $\widetilde P=\widetilde P^*\circ\cC^{-1}$ denote the induced observed laws. For every $t\in\cT$ and $a\in\cS_t$, by consistency and Equation~\eqref{eq:identical_marginals}
\begin{align*}
    P(T=t,Y=a)
    &= P^*(T=t,Y(t)=a) \\
    &= \sum_{\by\in\cS:y_t=a} P^*(T=t,\bYpot=\by) \\
    &= \sum_{\by\in\cS:y_t=a} \frac{1}{K}p_{\by} \\
    &= \frac{1}{K} \sum_{\by\in\cS:y_t=a}p_{\by} \\
    &= \frac{1}{K} (\bM_{\cS}\bp)_{t,a} \\
    &= \frac{1}{K} (\bM_{\cS}\widetilde{\bp})_{t,a} \\
    &= \frac{1}{K}\sum_{\by\in\cS:y_t=a}\widetilde p_{\by} \\
    &= \sum_{\by\in\cS:y_t=a}\frac{1}{K}\widetilde p_{\by} \\
    &= \sum_{\by\in\cS:y_t=a}\widetilde P^*(T=t,\bYpot=\by) \\
    &= \widetilde P^*(T=t,Y(t)=a) \\
    &= \widetilde P(T=t,Y=a).
\end{align*}
Hence $P=\widetilde P$ on the support of $\bYpot$. Since they are both zero outside the support, by identifiability in Definition~\ref{def:identifiability} it follows that $P^*=\widetilde P^*$ and, in particular, $\bp=\widetilde{\bp}$.

Therefore, for every $\bw \in \R^{\cS}$, $\bw^\top\bp=\bw^\top\widetilde{\bp}$. Since $\bp, \widetilde \bp \in \Delta(\cS)$ were arbitrary, Lemma~\ref{lem:marginalization}\ref{lem:marginalization_c} yields $\bw\in\mathrm{im}(\bM_{\cS}^\top)$. Furthermore, since $\bw\in\R^{\cS}$ was arbitrary, $\mathrm{im}(\bM_{\cS}^\top)=\R^{\cS}$. Thus $\bM_{\cS}$ has full column rank and, by Lemma~\ref{lem:marginalization}\ref{lem:marginalization_a}, $\bG_{\cS} = \bM_{\cS}^\top \bM_{\cS}$ is invertible.

Now suppose that $\bG_{\cS}$ is invertible. Fix $P^*\in\cPstar_{\cS}$ and define the probability vector
\[
    \bp = \lt(P^*(\bYpot=\by)\rt)_{\by\in\cS}.
\]
By the definition of the marginalization matrix $\bM_{\cS}$ in Lemma~\ref{lem:marginalization} and Assumptions~\ref{ass:consistency}, \ref{ass:ignorability}, and~\ref{ass:support},
\[
    \bq(\bp) := \bM_{\cS}\bp = \lt(P^*(Y(t)=a)\rt)_{t\in\cT,\,a\in\cS_t} = \lt(P(Y=a\mid T=t)\rt)_{t\in\cT,\,a\in\cS_t}.
\]
Since $\bG_{\cS}$ is invertible,
\[
    \bp = \bG_{\cS}^{-1}\bM_{\cS}^\top\bq(\bp). 
\]
It follows
\begin{align}
    P^*(\bYpot=\by)
    &= (\bp)_{\by} \nonumber\\
    &= \lt(\bG_{\cS}^{-1}\bM_{\cS}^\top
    \bq(\bp)\rt)_{\by} \nonumber\\
    &= \sum_{\by^\prime\in\cS}
    (\bG_{\cS}^{-1})_{\by,\by^\prime}
    \lt(\bM_{\cS}^\top\bq(\bp)\rt)_{\by^\prime} \nonumber\\
    &= \sum_{\by^\prime\in\cS}
    (\bG_{\cS}^{-1})_{\by,\by^\prime}
    \sum_{t\in\cT}\sum_{a\in\cS_t}
    (\bM_{\cS})_{(t,a),\by^\prime}
    \lt(\bq(\bp)\rt)_{t,a} \nonumber\\
    &= \sum_{\by^\prime\in\cS}
    (\bG_{\cS}^{-1})_{\by,\by^\prime}
    \sum_{t\in\cT}\sum_{a\in\cS_t}
    \indicator{y_t^\prime = a}
    P(Y=a\mid T=t) \nonumber\\
    &= \sum_{\by^\prime\in\cS}
    (\bG_{\cS}^{-1})_{\by,\by^\prime}
    \sum_{t\in\cT}
    P(Y=y_t^\prime\mid T=t).
    \label{eq:recovery_proof}
\end{align}
Since the right-hand side is a function only of the observed law $P$, the joint distribution $P^*(\bYpot)$ is identifiable. This also establishes the stated probability decomposition. By Assumption~\ref{ass:ignorability}, for every $t\in\cT$ and $\by\in\cS$,
\[
    P^*(T=t,\bYpot=\by)=P(T=t)P^*(\bYpot=\by).
\]
Since $P(T=t)$ is observed, the full-data law $P^*(T, \bYpot)$ is identifiable as well. This proves the other direction and thus the claimed equivalence. Finally, for any function $h:\cT\times\cY^{K}\to\R$, Equation~\eqref{eq:recovery_proof} yields
\begin{align*}
    \E_{P^*}[h(T,\bYpot)]
     &= \sum_{\by\in\cS}\sum_{t\in\cT} h(t,\by)P^*(T=t, \bYpot=\by) \\
    &= \sum_{\by\in\cS}\sum_{t\in\cT} h(t,\by)P(T=t)P^*(\bYpot=\by) \\
    &= \sum_{\by\in\cS}\sum_{t\in\cT} h(t,\by)P(T=t)
    \sum_{\by^\prime\in\cS}
    (\bG_{\cS}^{-1})_{\by,\by^\prime}
    \sum_{k\in\cT}P(Y=y_k^\prime\mid T=k) \\
    &= \sum_{\by,\by^\prime\in\cS}
    (\bG_{\cS}^{-1})_{\by,\by^\prime}
    \sum_{t\in\cT} h(t,\by)P(T=t)
    \sum_{k\in\cT}P(Y=y_k^\prime\mid T=k).
\end{align*}
This completes the proof.
\qed

\subsection{Proof of Corollary~\ref{cor:cardinality}}\label{app:cardinality}

Consider the marginalization matrix $\bM_{\cS}$ defined in Lemma~\ref{lem:marginalization}. By Lemma~\ref{lem:marginalization}\ref{lem:marginalization_b}, its kernel is given by
\begin{equation}\label{eq:cardinality_kernel}
    \mathrm{ker}(\bM_{\cS}) = \lt\{\bv \in \R^{\cS}: \sum_{\by \in \cS: y_t = a} v_{\by} = 0
    \ \text{for all } t \in \cT, a \in \cS_t \rt\}.
\end{equation}
For each $t \in \cT$, choose an arbitrary $b_t \in \cS_t$. We first show that $\mathrm{ker}(\bM_{\cS}) = \mathcal V_{\cS}$, where
\begin{equation}\label{eq:cardinality_reduced}
     \mathcal V_{\cS} := \lt\{\bv \in \R^{\cS}: \sum_{\by \in \cS} v_{\by} = 0,\ \sum_{\by \in \cS: y_t = a} v_{\by} = 0 \ \text{for all } t \in \cT,\ a \in \cS_t \setminus \{b_t\} \rt\}.
\end{equation}
From this, using identifiability of $P^*$ over $\cPstar_{\cS}$, we derive the claimed bounds.

Assume that $\bv \in \mathrm{ker}(\bM_{\cS})$ in  Equation~\eqref{eq:cardinality_kernel}. For any $t \in \cT$,
\[
    \sum_{\by \in \cS} v_{\by}
    = \sum_{a \in \cS_t} \sum_{\by \in \cS: y_t = a} v_{\by}
    = 0.
\]
This verifies the first condition of $\mathcal V_{\cS}$ in Equation~\eqref{eq:cardinality_reduced}. The second condition of $\mathcal V_{\cS}$ is satisfied trivially, thus $\bv \in \mathcal{V}_{\cS}$.

Conversely, let $\bv \in \mathcal{V}_{\cS}$. For every $t \in \cT$ and $a \in \cS_t\setminus \{b_t\}$ the condition in Equation~\eqref{eq:cardinality_kernel} holds. For the remaining $a = b_t$, we find
\begin{align*}
    \sum_{\by \in \cS: y_t = b_t} v_{\by}
    &= \sum_{\by \in \cS} v_{\by}
    - \sum_{a \in \cS_t \setminus \{b_t\}}
      \sum_{\by \in \cS: y_t = a} v_{\by} = 0,
\end{align*}
as both sums vanish by assumption. So $\bv \in \mathrm{ker}(\bM_{\cS})$ too. This proves that $\mathrm{ker}(\bM_{\cS}) = \mathcal V_{\cS}$.

The space $\mathcal{V}_{\cS}$ is defined through $|\cS|$ unknowns and $1+\sum_{t \in \cT}(|\cS_t|-1)$ linear equations. Since by assumption $P^*$ is identifiable over $\cPstar_{\cS}$, Theorem~\ref{thm:joint_identification} and Lemma~\ref{lem:marginalization}\ref{lem:marginalization_a} imply that $\mathcal{V}_{\cS} = \mathrm{ker}(\bM_{\cS}) = \{\bm{0}\}$. By \citep[Chapter 3, Corollary~$A^\prime$]{Lax2007}, we therefore must have
\[
    |\cS| \leq 1+\sum_{t \in \cT}(|\cS_t|-1),
\]
which proves Condition~\ref{cond:size_refined}. Since $\cS_t \subseteq \cY$ and $|\cY|=L$,
\[
    |\cS|
    \leq 1+\sum_{t \in \cT}(|\cS_t|-1)
    \leq 1+K(L-1),
\]
which proves Condition~\ref{cond:size_simple}.

Finally, we show that Condition~\ref{cond:size_refined} (and therefore Condition~\ref{cond:size_simple}) are not sufficient for identification. Let $K=3$, $L=2$, and
\[
    \cS = \{(0,0,0),(0,0,1),(1,1,0),(1,1,1)\}.
\]
Then $\cS_t=\{0,1\}$ for every $t \in \cT$,
\[
    |\cS| = 4 = 1+\sum_{t \in \cT}(|\cS_t|-1) = 1+K(L-1).
\]
However, $\bv=(1,-1,-1,1)^\top$ satisfies 
\[
    \sum_{\by \in \cS: y_t = 0} v_{\by} = 1-1 = 0,
    \qquad
    \sum_{\by \in \cS: y_t = 1} v_{\by} = -1+1 = 0,
\]
for every $t \in \cT$. Thus, by Equation~\eqref{eq:cardinality_kernel}, $\bv \in \mathrm{ker}(\bM_{\cS})$. However $\bv \neq \bm{0}$, therefore $\mathrm{ker}(\bM_{\cS}) \neq \{\bm{0}\}$ so $\bG_{\cS}$ is singular by Lemma~\ref{lem:marginalization}\ref{lem:marginalization_a} . By Theorem~\ref{thm:joint_identification}, $P^*$ is not identifiable over $\cPstar_{\cS}$.
\qed

\subsection{Proof of Theorem~\ref{thm:testable_implications}}\label{app:testable_implications}

Under Assumption~\ref{ass:support}, Theorem~\ref{thm:joint_identification} implies that
\[
    \pi_{\by}(P) = P^*(\bYpot = \by),
\]
for $\by \in \cS$. Hence $\pi_{\by}(P) \geq 0$. By Assumptions~\ref{ass:consistency} and~\ref{ass:ignorability},
\[
    \sum_{\by \in \cS: y_t = a} \pi_{\by}(P) = \sum_{\by \in \cS: y_t = a} P^*(\bYpot = \by) =  P^*(Y(t) = a)
    = P(Y = a \mid T = t),
\]
for every $t \in \cT$ and $a \in \cY$. Hence, Equation~\eqref{eq:testable_conditions} is satisfied.
Consequently, failure of Equation~\eqref{eq:testable_conditions} falsifies Assumption~\ref{ass:support} under consistency and strong ignorability.

To establish sharpness, let $P$ satisfy Equation~\eqref{eq:testable_conditions} and overlap, $P(T= t) > 0$ for all $t \in \cT$. We will show that there exists a $P^* \in \cPstar_{\cS}$ for which $P = P^* \circ \cC^{-1}$.

First, by Equation~\eqref{eq:testable_conditions}, summing the marginal probabilities for any fixed $t \in \cT$
\[
    \sum_{\by \in \cS} \pi_{\by}(P) = \sum_{a \in \cS_t} \sum_{\by \in \cS: y_t = a} \pi_{\by}(P) =  \sum_{a \in \cS_t} P(Y = a \mid T = t) = 1.
\]
As $\pi_{\by}(P) \geq 0$ for all $\by \in \cS$, we define the full data law
\[
    P^*(T = t, \bYpot = \by)
    = P(T = t)\pi_{\by}(P),
\]
for all $t \in \cT, \by \in \cS$, with zero mass outside. Under $P^*$, $T$ is independent of $\bYpot$. Overlap holds by assumption and $\bYpot \in \cS$ by construction. Therefore, $P^* \in \cPstar_{\cS}$. Moreover,
\begin{align*}
    P^*(T = t, Y(t) = a)
    = P(T = t)\sum_{\by \in \cS: y_t = a} \pi_{\by}(P) = P(T = t, Y = a),
\end{align*}
for every $t \in \cT$ and $a \in \cY$. Hence $P^* \circ \cC^{-1} = P$. This completes the proof.

\subsection{A Technical Lemma}\label{app:marginalization}


The following lemma introduces a zero--one matrix that maps joint probabilities to their marginals, connecting the properties of this matrix to our identification results. A similar construction appears in the instrumental variable model of \citet{heckman2018unordered} to describe treatment choices. Here, we use the Gram representation $\bG_{\cS}=\bM_{\cS}^{\top}\bM_{\cS}$ to obtain an interpretable identification condition based on Hamming distances between principal strata and derive closed-form identification formulas.

\begin{lemma}\label{lem:marginalization}
Let $\varnothing \ne \cS \subseteq \cY^{K}$ be finite and write $\cS_t = \{y_t: \by \in \cS\}$. Let $\Delta(\cS) = \{\bp \in [0, 1]^{\cS}: \sum_{\by \in \cS}p_{\by} = 1\}$ be the probability simplex of $\bYpot$ supported on $\cS$. Define the zero-one marginalization matrix $\bM_{\cS} \in \{0, 1\}^{(\sum_{t \in \cT} |\cS_t|) \times |\cS|}$ by
\begin{equation*}
    (\bM_{\cS})_{(t,a),\by} := \indicator{y_t = a},
\end{equation*}
for $t \in \cT$, $a \in \cS_t$, and $\by \in \cS$.
The matrix $\bM_{\cS}$ maps a joint distribution $\bp \in \Delta(\cS)$ to the vector of its marginal distributions $\bq(\bp)$, i.e.,
\[
    \bq(\bp) = \bM_{\cS} \bp = \lt(\sum_{\by \in \cS: y_t = a} p_{\by} \rt)_{t \in \cT, a \in \cS_t}.
\]
Then the following statements hold:
\begin{enumerate}[label=\textnormal{(\alph*)}, ref=\textnormal{(\alph*)}]
    \item \label{lem:marginalization_a} $\bG_{\cS}=\bM_{\cS}^{\top}\bM_{\cS}$. In particular, $\bG_{\cS}$ is invertible if and only if $\bM_{\cS}$ has full column rank $\mathrm{rank}(\bM_{\cS})=|\cS|$, or equivalently, $\mathrm{ker}(\bM_{\cS})=\{\bm{0}\}$.
    \item \label{lem:marginalization_b} The kernel of $\bM_{\cS}$ is given by
    \[
        \mathrm{ker}(\bM_{\cS}) = \lt\{\bv \in \R^{\cS}: \sum_{\by \in \cS: y_t = a} v_{\by} = 0
        \ \text{for all } t \in \cT, a \in \cS_t \rt\}.
    \]
    Consequently, $\mathbf{1}^\top\bv = 0$ for all $\bv \in \mathrm{ker}(\bM_{\cS})$.
    \item \label{lem:marginalization_c} Fix $\bw \in \R^{\cS}$. If, for all $\bp, \bp^\prime \in \Delta(\cS)$ with $\bM_{\cS}\bp=\bM_{\cS}\bp^\prime$ we have $\bw^\top \bp = \bw^\top \bp^\prime$, then $\bw\in\mathrm{im}(\bM_{\cS}^\top)$.
\end{enumerate}
\end{lemma}

\begin{proof}
We first prove (a). We show that $\bG_{\cS} = \bM_{\cS}^\top \bM_{\cS}$. Indeed, for any $\by,\by^\prime \in \cS$,
\begin{align*}
    (\bM_{\cS}^\top\bM_{\cS})_{\by,\by^\prime}
    = \sum_{t \in \cT} \sum_{a \in \cS_t}
    \indicator{y_t = a}\indicator{y_t^\prime = a} = \sum_{t \in \cT}\indicator{y_t = y_t^\prime}
    = (\bG_{\cS})_{\by,\by^\prime}.
\end{align*}

Next, we show that $\bG_{\cS}$ is invertible if and only if $\bM_{\cS}$ has full column rank. For every $\bv \in \R^{\cS}$,
\[
    \bv^\top \bG_{\cS} \bv = \bv^\top \bM_{\cS}^\top \bM_{\cS} \bv = \|\bM_{\cS}\bv\|_2^2.
\]
Therefore $\bv \in \mathrm{ker}(\bG_{\cS})$ if and only if $\bv \in \mathrm{ker}(\bM_{\cS})$. Consequently, 
\begin{align*}
    \bG_{\cS} \text{ is invertible}
    &\iff \mathrm{ker}(\bG_{\cS}) = \{\bm{0}\} \\
    &\iff \mathrm{ker}(\bM_{\cS}) = \{\bm{0}\} \\
    &\iff \mathrm{rank}(\bM_{\cS}) = |\cS|,
\end{align*}
where the last equivalence follows from theorem \citep[Chapter~3, Theorem 2]{Lax2007}.

Next, we prove (b). Let $\bv = (v_{\by})_{\by \in \cS} \in \R^{\cS}$, then
\[
    (\bM_{\cS} \bv)_{t, a} = \sum_{\by\in \cS} \indicator{y_t = a} v_{\by} = \sum_{\by \in \cS: y_t = a} v_{\by}.
\]
Therefore
\[
    \bv \in \mathrm{ker}(\bM_{\cS}) \iff \bM_{\cS} \bv = 0 \iff \sum_{\by \in \cS: y_t = a} v_{\by} = 0,
\]
for all $t \in \cT, a \in \cS_t$, which proves the kernel characterization.
Consequently, if $\bv \in \mathrm{ker}(\bM_{\cS})$, then, for any $t \in \cT$,
\begin{align*}
    \mathbf{1}^\top\bv = \sum_{\by \in \cS} v_{\by} = \sum_{a \in \cS_t}  \sum_{\by \in \cS: y_t = a} v_{\by} = 0.
\end{align*}

Finally, we prove (c). Fix $\bw \in \R^{\cS}$ and assume that $\bw \notin \mathrm{im}(\bM_{\cS}^\top)$. By standard linear algebra \citep[Chapter~3, Theorem~$5^\prime$]{Lax2007}, $\mathrm{im}(\bM_{\cS}^\top) = \ker(\bM_{\cS})^\perp$. That is, the image of $\bM_{\cS}^\top$ is the orthogonal complement of the kernel of $\bM_{\cS}$. This implies that there exists $\bv = (v_{\by})_{\by \in \cS} \in \ker(\bM_{\cS}) \setminus \{\bm{0}\}: \bw^{\top}\bv\neq 0$. Let $\bp = (p_{\by})_{\by \in \cS} \in \Delta(\cS)$ have strictly positive entries, for example $p_{\by}=1/|\cS|$ for every $\by \in \cS$. Define
\begin{align*}
    \bp^\prime = \bp + \alpha \bv, \quad \text{with} \ 0 < \alpha < \frac{\min_{\by \in \cS} p_{\by}}{\max_{\by \in \cS} \abs{v_{\by}}}.
\end{align*}
Note that such an $\alpha$ exists. Indeed, by assumption, $0 < p_{\by}$ for all entries and $\bv$ has at least one non-zero entry as $\bv\neq \bm{0}$. Observe,
\begin{align*}
    \min_{\by \in \cS} p^\prime_{\by} \geq \min_{\by \in \cS} p_{\by} - \alpha \max_{\by \in \cS} \abs{v_{\by}} > 0,
\end{align*}
and,
\begin{align*}
    \sum_{\by \in \cS} p^\prime_{\by} = \sum_{\by \in \cS} p_{\by} + \alpha \sum_{\by \in \cS} v_{\by} = \sum_{\by \in \cS} p_{\by} = 1,
\end{align*}
where the second equality follows from $\mathbf{1}^\top\bv=0$ by $\bv\in\mathrm{ker}(\bM_{\cS})$ and part (b). Thus, $\bp^\prime\in\Delta(\cS)$.

Because $\bv \in \text{ker}(\bM_{\cS})$,
\begin{align*}
    \bM_{\cS} \bp^\prime = \bM_{\cS}(\bp + \alpha \bv) = \bM_{\cS} \bp.
\end{align*}
By (c), this implies that $\bw^\top\bp=\bw^\top\bp^\prime$. Consequently, 
\begin{align*}
     0 = \bw^\top \bp^\prime - \bw^\top \bp
    = \alpha \bw^\top \bv \neq 0,
\end{align*}
as $\alpha \neq 0$ by construction and $\bw^\top \bv \neq 0$ by the choice of $\bv$, since $\bw \notin \mathrm{im}(\bM_{\cS}^\top) = \ker(\bM_{\cS})^\perp$. This yields a contradiction and proves (c).

\end{proof}

\end{document}